\documentstyle [11pt]{article}
\newtheorem{theorem}{Theorem}[section]

\newtheorem{prop}[theorem]{Theorem}

\newtheorem{coroll}[theorem]{Corollary}

\newtheorem{lemm}[theorem]{Lemma}

\newtheorem{anoprop}[theorem]{Proposition}

\newenvironment{mproof}{\begin{trivlist}\item[]{\em
Proof: }}{\hfill$\Box$\end{trivlist}}

\newtheorem{eg}{\rm\sl \uppercase{Example}}[section]

\def \IR{\hbox{{\rm I}\kern-.2em\hbox{{\rm R}}}}
\def \iR{\hbox{{\sevenrm I\kern-.2em\hbox{\sevenrm R}}}}
\def \IN{\hbox{{\rm I}\kern-.2em\hbox{\rm N}}}

\def \IC{\hbox{{\rm I}\kern-.6em\hbox{\bf C}}}
\def \IQ{\hbox{{\rm I}\kern-.6em\hbox{\bf Q}}}
\def \ZZ{\hbox{{\rm Z}\kern-.4em\hbox{\rm Z}}}

\newcommand{\bmi}{\mbox{\boldmath $i$}}
\newcommand{\bmj}{\mbox{\boldmath $j$}}

\begin{document}

\begin{center}

{\Large\bf Infinite Lexicographic Products of Triangular Algebras}

\end{center}
\vspace{.3in}
\begin{center}
{\Large S.C. Power}

\vspace{.3in}
\it Department of Mathematics\\
Lancaster University\\
England LA1 4YF
\rm

\vspace{.3in}

ABSTRACT
\end{center}
{\small Some new connections are given between
linear orderings and
triangular operator algebras. A lexicograhic product is defined
for triangular operator algebras
and the Jacobson radical of an infinite lexicographic product
of upper triangular matrix algebras is determined.}

\vspace{.5in}
\begin{center}
1. Introduction
\end{center}

For certain triangular operator algebras we define a
natural lexicographic product $A_1 \star A_2$, which is particularly
pertinent to the category of maximal triangular operator
algebras, and we discuss the associated infinite lexicographic
products over linear orderings. Our main purpose is the identification
of the Jacobson radical of the lexicographic products
\[
A(\Omega, \nu) = \prod_{w \in \Omega } \star \ T_{n_\omega}
\]
where $\Omega$ is a countable linear
ordering,  where $\nu : \omega \to n_\omega$
is a map from  $\Omega$ to $\{2, 3, 4, \dots\}$, and
where the factors $T_n$ are
complex upper triangular matrix algebras.

The algebra $A(\Omega, \nu)$
is the closure in a natural operator norm of an algebraic direct limit
$A_\infty(\Omega, \nu)$ of upper triangular matrix algebras.  Let us say
that a complex algebra $A$ has an elementary radical decomposition if,  like
$T_n$, $A$ can be written as a direct sum $A = C + rad (A)$ where $C$ is
abelian and where $rad  (A)$ is the Jacobson  radical. Although
$A_\infty(\Omega, \nu)$ admits such a decomposition we show that
$A(\Omega, \nu)$ has an elementary radical decomposition if and only if
$\Omega$ is well-ordered.  Furthermore, in the case of a general
linear ordering we determine the
Jacobson radical and its semisimple quotient.  If $\Omega = \ZZ_+, \ZZ_-$ or
$\ZZ$, then we recover the refinement limit algebras, the standard limit
algebras and the (proper) alternation algebras,
respectively.  For discussion of these
fundamental examples and related matters see \cite{scp-book}.
The principal
result we need to invoke is Donsig's recent characterisation of semisimple
triangular limit algebras.  For this see Donsig \cite{don} or \cite{scp-book}.

The algebras $A_\infty(\Omega, \nu)$ and
$A(\Omega, \nu)$ provide many new isomorphism types of (strongly
maximal) triangular algebras and they exhibit new phenomena.
In particular, in the $2^\infty$ UHF C*-algebra
there are uncountably many lexicographic
products. Also
there are unital limit algebras of the
form ${\displaystyle \lim_\to T_{2^k}}$ with non-abelian outer
automorphism group modulo pointwise inner automorphisms. A detailed
account of  classification, isomorphisms and automorphisms will be
given elsewhere. Significantly, it is possible to reduce many such
considerations to the well-developed theory of linear orders and
their automorphisms. (See Rosenstein  \cite{ros}.) We remark that
other connections between linear orderings
and triangular operator algebras have appeared in \cite{muh-saito-sol-1}
and \cite{orr}

\begin{center}
2. Lexicographic Products
\end{center}

Let $\Omega, \nu$ be as above, let $F \subseteq \Omega$ be a finite
subset, say $w_1 < w_2 < \dots < w_k$, and let $w_t < w < w_{t+1}$, for some
$t$. Set $G = F \cup \{w\}$, $n_F = n_{w_1} n_{w_2}\dots n_{w_k}$, and $n_G =
n_\omega n_F$.   Define a unital algebra injection $\phi_{F,G} : T_{n_F} \to
T_{n_G}$ as follows. View $T_{n_F}$ as the (maximal triangular) subalgebra
of $M_{n_1} \otimes \dots \otimes M_{n_{w_k}}$ which is spanned by the
matrix units
\[
e_{{\bmi},{\bmj}} = e_{i_1, j_1} \otimes \dots \otimes e_{i_k,j_k}
\]
where the multi-index ${\bmi} = (i_1, \dots, i_k)$
precedes ${\bmj} = (j_1, \dots,
j_k)$ in the lexicographic ordering. Thus either ${\bmi} =  {\bmj}$ or the
first
$i_p$ differing from $j_p$ is strictly less than $j_p$. Similarly  identify
$T_{n_G}$ for the ordered subset $G$ and  set $\phi_{F,G}$ to be the linear
extension of the correspondence
\[
e_{{\bmi},{\bmj}} \to \sum^{n_\omega}_{s=1} e_{i_1, j_1} \otimes \dots \otimes
e_{i_t, j_t} \otimes e_{s,s} \otimes e_{i_{t-1},j_{t+1}} \otimes \dots \otimes
e_{i_{k},j_{k}}
\]

In a similar way (or by composing maps of the above type) define
$\phi_{F,G}$ for $F \subseteq G$, general finite subsets. These maps are
isometric and so determine the Banach  algebra
\[
A(\Omega, \nu) = {\lim_{\to}}_{F \in {\cal F}} T_{n_F}
\]
where the direct limit is taken over the directed set ${\cal F}$ of finite
subsets of $\Omega$. Since each $\phi_{F,G}$ has an extension to a
C*-algebra injection from $M_{n_F}$ to $M_{n_G}$ it follows that we may
view $A(\Omega, \nu)$ as a closed unital subalgebra of the UHF
C*-algebra $B(\Omega, \nu) = {\displaystyle \lim_\to M_{n_F}}$.

Let $A_1 \subseteq M_n, A_2 \subseteq M_m$ be triangular digraph
algebras, and let $A^0_1$ be the maximal ideal of $A_1$ which is disjoint
from the diagonal masa $A_1 \cap A_1^*$. Define the lexicographic
product $A_1 \star A_2$ to be the triangular digraph subalgebra of $M_n
\otimes M_m$ given by
\[
A_1 \star A_2 = (A_1 \cap A_1^*) \otimes A_2 + A^0_1 \otimes C^\ast (A_2).
\]
Note that there are natural inclusions $A_1 \to A_1 \star A_2$, $A_2 \to A_1
\star A_2$. Furthermore $\star$ is associative, and the injections
$\phi_{F,G}$ of the last paragraph can be identified as the natural inclusion
map
\[
\prod_{w\in F} \star \ T_{n_\omega} \to \prod_{w \in G} \star \ T_{n_\omega}
\]
for the finite lexicographic products for $F$ and $G$. Although $T_n \star T_m
= T_{nm} = T_m \star T_n$, the operation $\star$ is not commutative in general.

We can use the same formula as above to define a general lexicographic
product $A_1 \star A_2$ whenever $A_1$ is an operator algebra admitting a
decomposition $A_1 = A_1 \cap A_1^\ast + A^0_1$ where $A_1 \cap A^*_1$ is a
maximal abelian self-adjoint subalgebra of $A_1$ and $A^0_1$ is an ideal
which is the kernel of a contractive homomorphism $A_1 \to A_1 \cap
A^\ast_1$.  For definiteness, we take  $A_1 \star A_2$ to be the normed
subalgebra of the injective C*-algebra tensor product.
The algebras  $A(\Omega, \nu)$ themselves fall into this
category, with $A(\Omega, \nu)^0$ equal to the ideal ${\displaystyle \lim_\to
T_{n_F}^0}$. Clearly
\[
A(\Omega, \nu) \star A(\Lambda, \mu) = A(\Omega + \Lambda, \nu + \mu)
\]
where $\Omega + \Lambda$ is the order sum of $\Omega$ and $\Lambda$.

As with the cases $\Omega = \ZZ, \ \ZZ_+, $ and \ $\ZZ_-$,
the Gelfand space of $C = A
\cap A^\ast$ is naturally identifiable with the Cantor space $ X
= \prod_{\omega \in \Omega} [n_\omega]$ where $[n_\omega] =
\{1,...,n_\omega \}$. Write $x \sim y$ if $x = (x_\omega)$ and $y = (y_\omega)$
are points in $X$ with $x_\omega = y_\omega$ for all but
finitely many $\omega$.  This
equivalence relation carries a natural topology and $B(\Omega, \nu)$ can be
viewed as the groupoid C*-algebra of this topological
relation. (See \cite{ren}). From this perspective
$A(\Omega, \nu)$ is the triangular subalgebra determined by the
lexicographic subrelation.
Similarly if  $A(G_\omega) \subseteq T_{n_\omega}$ are
unital digraph algebras then the lexicographic product
$A = \prod_{w \in \Omega} \star A(G_\omega)$
can also be viewed as a semigroupoid
algebra associated with the natural subrelation of $\sim$.

\begin{center}
3. The Jacobson Radical
\end{center}

\noindent {\bf Lemma 1.}\ \
\it
Let $\Omega$ be well-ordered. Then  $rad(A(\Omega, \nu))$ coincides with the
maximal diagonal disjoint ideal  ${\displaystyle A(\Omega, \nu)^0
= {\lim_{\to}}_{F\in{\cal F}} T_{n_F}^0}$.  Furthermore, if $B$ is an AF
C*-algebra
then the Jacobson radical of $A(\Omega, \nu) \otimes B$ is $A (\Omega,
\nu)^0 \otimes B$.
\rm

\begin{mproof}
Identify the building block algebras $T_{n_\omega}$ with their images in $A =
A(\Omega,
\nu)$. Assume that the conclusion is false and let $w_1 \in \Omega$ be the
least element $w$ for which  $T^0_{n_\omega}$ is not contained in $rad A$. Let
$\Omega_1 = \{ w \in \Omega : w < w_1\}, \ \Omega_2 = \Omega \setminus
\ \Omega_1, \ A_i = A(\Omega_1, \nu), \ $ for $ i = 1, 2,$
and identify $A$ with
$ A_1 \star A_2 = C_1 \otimes A_2 + A_1^0 \otimes C^\ast(A_1)$, where  $C_1 =
A_1 \cap A^\ast_1$, and $A^0_1$ is the maximal diagonal disjoint ideal of
$A_1$. From the definition of $w_1$ it follows that $A^0_1 \otimes \ \IC
\subseteq rad A$, and  hence that $A^0_1 \otimes C^\ast (A_1) \subseteq rad
A$, so that
\[
A/rad A = (C_1 \otimes A_2) /rad A.
\]
Write $A_2$ as $T_{n_{w_1}} \star A(\Omega_3, \nu)$ where $\Omega_3 =
\Omega_2 \setminus \{w_1\}$ and observe that a matrix unit $a =
e_{i,j}$ in $T_{n_{w_1}}$, with $i < j$, satisfies $(ab)^p = 0$ for all $b$ in
$A_2$ where $p = n_{\omega}$.  Indeed $A_2$ can be
viewed as a  $p \times p$ block upper
triangular algebra, with  the element $a$ appearing in the
$i,j$ position. It follows that $(ab)^p = 0$ for all $b$ in $C_1 \otimes
A_2$, and hence that $a + rad A$ belongs to the radical of the quotient. Since
the quotient is semisimple, $a$ belongs to $ rad A$.
Thus $T_{n_{w_1}}^0 \subseteq rad A $, a contradiction.

For the final part note that the proof above can be repeated, starting with the
least $w_1$ for which $T_{n_{w_1} }\otimes B$ is not contained in the radical.
\end{mproof}

\noindent {\bf Lemma 2.}\ \
\it
\ Suppose that $\Omega$ does not have a first element. Then
the algebra $A(\Omega, \nu)$ is semisimple.
\rm

\begin{mproof}
By Donsig's theorem it will be enough to show that for each
matrix unit $e$ in $T_{n_\omega} \subseteq A$ there is a link for $e$, that is,
there is a
matrix unit $f$ such that $f^\ast f \leq ee^\ast$ and $ff^\ast \leq e^\ast e$.
By
the hypothesis, for any such element $e$ there exists $w_1 < w$, and with
respect to the identification of $T_{n_{w_1}} \star T_{n_\omega}$ in
$M_{n_{w_1}}
\otimes M_{n_\omega}$ the matrix unit $e$ is
identified with $I \otimes e$. Let $f$ be the matrix
unit in $T_{n_{w_1}} \star T_{n_\omega}$ with initial projection  $e_{22}
\otimes
ee^\ast$ and final projection $e_{11} \otimes e^\ast e$. Then $f$ is a link for
$e$.
\end{mproof}

\noindent {\bf Theorem 3.}\
\it
\ Let $A = A(\Omega, \nu), $
and for $i = 1,2,$ let $A_i = A(\Omega_i, \nu_i)$,
where $\Omega  = \Omega_1 + \Omega_2$ is the order sum decomposition
for which $\Omega_1$ is the maximal well-ordered initial segment and
$\nu_i = \nu|\Omega_i, i = 1, 2$. Then the Jacobson radical of $A$ is the ideal
$A_1^0 \otimes C^\ast(A_2)$ and $A/rad(A)$ is completely isometrically
isomorphic to $(A_1 \cap A^\ast_1)  \otimes A_2$.
\rm

\begin{mproof}
Let $J$  be the ideal $A^0_1 \otimes C^\ast(A_2)$ of  $A = A_1 \star A_2$,
where $A$ is identified, as usual, with a subalgebra of $A_1 \otimes C^\ast
(A_2)$. Since $J$ is also an ideal of $A_1 \otimes C^\ast(A_2)$ it follows that
the natural  map $ A/J \to (A_1 \otimes C^\ast(A))/J$ is a complete
isometry, and so we have a completely isometric injection
\[
\alpha : A/J \to
(A_1/A^0_1) \otimes C^\ast(A) = (A_1 \cap A^\ast_1) \otimes C^\ast(A_2).
\]
But $A_1 \star A_2 = (A_1 \cap A^\ast_1) \otimes A_2 + A^0_1 \otimes C^\ast
(A_2)$ and so it follows that the range of $\alpha$ is precisely  $(A_1 \cap
A^\ast_1) \otimes A_2$.

By Lemma 2 $A_2$ is semisimple and so $(A_1 \cap A^\ast_1) \otimes A_2$ is
semisimple. Thus $rad(A)$ is contained in $J$. But on the other hand, by
Lemma 1, $J$ is precisely the radical of $A_1 \otimes C^\ast(A_2)$ and so it
follows that $J$ is contained in the radical of the subalgebra $A$.
\end{mproof}

  The arguments above are much more
generally applicable. For example if
${
A = \prod_{w\in\Omega}}\star A(G_\omega)$
where each  $A(G_\omega)$ is a triangular digraph
algebra with a proper elementary radical decomposition, then the
Jacobson radical is similarly identified.

  The algebras of Theorem 3 provide many specific examples
of maximal triangular subalgebras of groupoid C*-algebras in the sense
of Muhly and Solel \cite{muh-sol}. On the other hand
quite different classes of maximal
triangular subalgebras of C*-algebras arise from
lexicographic products involving
function algebras such as the disc algebra $A(D)$.
In this regard note that $A(D) \star T_n$ is semisimple for any $n$, and
in fact it can be shown that
all lexicographic products of such algebras are
semisimple.

\begin{center}
4. Classification
\end{center}

If $\gamma : \IQ \to \IQ$ is an order bijection, and $A(\IQ , 2^\infty)$ is the
lexicographic product with $n_q = 2$ for all $q$, then there is a natural
induced isometric automorphism $\alpha_\gamma$ of $A(\IQ , 2^\infty)$. In
fact with respect to the natural inclusion
\[
A(\IQ, 2^\infty) \to \prod_{q \ \in \ \IQ} \otimes M_2(\IC)
\]
the map $\alpha_\gamma$ is simply the restriction of the shift automorphism of
$B(\Omega, 2^\infty)$ induced by $\gamma$.  It can be shown that
modulo the approximately inner automorphisms all isometric automorphisms
arise this way. Thus $Out_{isom} (A(\IQ , 2^\infty)) \simeq Aut(\IQ)$. This
seems to be the first example of a unital limit of upper triangular matrices
for
which the outer isometric automorphism group
is not abelian (cf. \cite{scp-outer}).
More generally it can be shown,
by an analysis of
the closed semi-orbits
for the associated lexicographic semigroupoid, that $A(\IQ, \nu)$ and $A(\IQ,
\mu)$ are isomorphic algebras if and only if there is a bijection $\gamma$
with $n_{\mu(q)} = n_{\nu(\gamma(q))}$ for all $q$ in $\IQ$.
Furthermore there is a
complete classification of all  the algebras $A(\Omega, \nu)$ which can be
made by a similar analysis.  Associate with each $\omega \in \Omega$ the
maximal order interval $I_\omega$, containing $\omega$,
such that $I_\omega$ is isomorphic to
$\ZZ, \ \ZZ_+$ or \ $\ZZ_-$, or is finite.
Each such interval
has associated with
it an upper triangular
matrix algebra or an alternation algebra for the data $n_\omega$, for $\omega
\in I_\omega$. Two lexicographic product algebras are isomorphic if and only
if,
modulo a parameter change $\gamma$, the linearly
ordered sets of maximal intervals
are isomorphic, and the associated alternation algebras agree. Thus the
only essential variation in presentation of an isomorphism class is that
which is already present in the case of alternation algebras. See
\cite{hop-scp},
\cite{poon}, \cite{scp-outer}, \cite{scp-lex}.

\end{document}